\documentclass{amsart} 
\usepackage{amssymb,amscd,verbatim,latexsym, graphicx}
\usepackage{enumerate}
\usepackage{amsmath,amsthm}
\usepackage{float}
\usepackage{color}
\newcommand{\CC}{\mathbb C} 
 
 \newcommand{\RR}{\mathbb R}
 
\newcommand{\ZZ}{\mathbb Z}

\newcommand{\CI}{{\mathcal C}^{\infty}}
\newcommand{\CIc}{{\mathcal C}^{\infty}_{\text{c}}}
\newcommand\pa{{\partial}}

\newcommand{\loc}{\operatorname{loc}}

\newcommand{\maC}{\mathcal C}
\newcommand{\maD}{\mathcal D}

\newcommand{\maK}{\mathcal K}

\newcommand{\maV}{\mathcal V}

\newcommand{\maT}{\mathcal T}

\newcommand\ttra{tetrahedralization}

\newcommand{\ie}{{\em i.e., }}

\newtheorem{theorem}{Theorem}[section]
\newtheorem{proposition}[theorem]{Proposition}

\theoremstyle{definition}

\theoremstyle{remark}

\numberwithin{figure}{section}
\numberwithin{table}{section}

\newcommand\Kond[2]{{\mathcal K}^{#1}_{#2}}
\newcommand{\psO}{\partial_{\operatorname{sing}} \Omega}
\newcommand{\BK}{Babu\v{s}ka-Kondratiev}

\definecolor{DarkRed}{rgb}{0.8,0.3,0.6}  

\definecolor{MidnightBlue}{rgb}{.1,.1,0.6}

\definecolor{Green}{rgb}{0.2,0.7,0.3}

\title[FEM on 3D polyhedral domains]{Anisotropic
  regularity and optimal rates of convergence for the Finite Element
  Method on three dimensional polyhedral domains}

\author[C. B\u{a}cu\c{t}\u{a}]{Constantin B\u{a}cu\c{t}\u{a}}
\address{University of Delaware,
        Department of Mathematical Sciences,
        501 Ewing Hall,
        Newark, DE 19716-2553, USA.}
\email{bacuta@math.udel.edu}

\author[A. Mazzucato]{Anna L. Mazzucato}
\address{Pennsylvania State University,
        Department of Mathematics,
        University Park, PA 16802, USA.}
\email{alm24@psu.edu}

\author[V. Nistor]{Victor Nistor}
\address{Pennsylvania State University,
        Department of Mathematics,
        University Park, PA 16802, USA.}
\email{nistor@math.psu.edu}

\thanks{The work of C. Bacuta is partially supported by NSF
  DMS-0713125.  The work of A. Mazzucato is partially supported by NSF
  grant DMS 1009713, 1009714, and by a Michler Fellow appointment at
  Cornell University.  The work of V. Nistor is partially supported by
  NSF grant DMS 1016556.}

\begin{document}

\begin{abstract}
We consider the model Poisson problem $-\Delta u = f \in \Omega$, $u =
g$ on $\pa \Omega$, where $ \Omega $ is a bounded polyhedral domain in
$\RR^n$. The objective of the paper is twofold. The first objective is
to review the well posedness and the regularity
of our model problem using appropriate weighted spaces for the data
and the solution. We use these results to derive the domain of
  the Laplace operator with zero boundary conditions on a concave
  domain, which seems not to have been fully investigated before. We
  also mention some extensions of our results to interface problems
  for the Elasticity equation. The second objective is to illustrate
how anisotropic weighted regularity results for the Laplace operator
in 3D are used in designing efficient finite element discretizations
of elliptic boundary value problems, with the focus on the efficient
discretization of the Poisson problem on polyhedral domains in
$\RR^3$, following {\em Numer. Funct. Anal. Optim.}, 28(7-8):775--824, 2007. The anisotropic weighted regularity
results described and used in the second part of the paper are a
consequence of the well-posedness results in (isotropically) weighted
Sobolev spaces described in the first part of the paper. The
  paper is based on the talk by the last named author at the Congress
  of Romanian Mathematicians, Brasov 2011, and is largely a survey
  paper.
\end{abstract}

\maketitle

%\tableofcontents

\section*{Introduction}

Let $\Omega \subset \RR^n$ be an open, bounded set. Consider the
boundary value problem
\begin{equation}\label{eq.BVP}
  \begin{cases} \; \Delta u = f & \text{ in }\Omega  
      \\ u\vert_{\pa \Omega} = g, & \text{ on } \Omega,
\end{cases}
\end{equation}
defined on a bounded domain $\Omega \subset \RR^n$, where $\Delta$ is
the Laplacian $\Delta = \sum_{i=1}^d \pa_i^2$. When $\pa \Omega$ is
smooth, it is well known that this Poisson problem has a unique
solution $u \in H^{m+1}(\Omega)$ for any $f \in H^{m-1}(\Omega)$ and
$g \in H^{m+1/2}(\pa \Omega)$ \cite{Evans, LionsMagenes1, Taylor1}.
Moreover, $u$ depends continuously on $f$ and $g$. This result is the
{\em classical well-posedness of the Poisson problem} on smooth
domains.

On the other hand, when $\Omega$ is not smooth, it is also well known
\cite{CDSchwab, Dauge, JerisonKenig, Kondratiev67, KMRossmann} that there
exists $s = s_{\Omega}$ such that $u \in H^{s}(\Omega)$ for any $s <
s_{\Omega}$, but $u \not\in H^{s_{\Omega}}(\Omega)$ in general, even
if $f$ and $g$ are smooth functions defined in a neighborhood of
$\Omega$. For instance, if $\Omega$ is a polygonal domain in two
dimensions, then $s_{\Omega} = 1+\pi/\alpha_{MAX}$, where
$\alpha_{MAX}$ is the largest interior angle of $\Omega$
\cite{Kondratiev67}. See also Wahlbin's paper \cite{Wahlbin84}.

In view of applications to the Finite Element Method, we restrict out
attention to domains {\em with polyhedral structure}.  These are
natural non-convex generalizations of classical $n$-dimensional
polyhedra that allow for curved boundaries, cracks (i.e., internal
faces), and non-smooth boundary points touching a smooth part of the
boundary. We refer to \cite{BMNZ} for a precise formulation.

There exist a very large number of papers devoted to boundary
value problems on non-smooth domains. While it is impossible to
mention them all, let us at least mention the papers of Arnold, Scott,
and Vogelius \cite{ASV}, Babuska and Guo \cite{BabuGuo2},
B\u{a}cu\c{t}\u{a}, Bramble, and Xu \cite{bacutaBX}, Jerison and Kenig
\cite{JerisonKenig}, Kondratiev \cite{Kondratiev67}, Kozlov, Mazya,
and Rossmann \cite{KMRossmann}, Mitrea and Taylor \cite{MMTaylor},
Rossmann \cite{Rossmann}, Verchota \cite{Verchota}, and many
others. Other results specific to numerical methods for polyhedral
domains are contained in the papers of Apel and Dobrowolski
\cite{ApelDobro}, Costabel, Dauge, and Nicaise \cite{CDN}, Costabel,
Dauge, and Schwab \cite{CDSchwab}, Dauge \cite{Dauge}, Demkowicz,
Monk, Schwab, and Vardapetyan \cite{Monk00}, Elschner
\cite{Elschner1}, Guo and Schwab \cite{GuoSchwab}, and many
others. Further results and references can be found in the
aforementioned papers, as well as in the the monographs of Grisvard
\cite{Grisvard2} as well as the recent book \cite{NP}. Regularity for
polyhedral domains is useful in designing fast solvers for numerical
methods \cite{ApelMG, BrannickHengguang}. See also \cite{Apel99,
  ApelDobro, BBP3, BlumDobro, BrennerSung, LeeHengguang, HNS, Felli1,
  Felli2, Hengg2D, Manteuffel2006} for more applications of these
techniques to other types of Partial Differential Equations and
numerical methods.

In this paper, we shall review the results from \cite{BMNZ, 3D1, 3D2,
  HMN}, and \cite{MNelast}, which make use of the natural stratified
space structure on $\Omega$. This leads, by successive conformal
changes of the metric, to a metric for which the smooth part of
$\overline{\Omega}$ becomes a smooth manifold with boundary whose
double is complete. The resulting Sobolev spaces defined by the new
metric will lead to spaces on which the Poisson problem is well-posed.

We restrict for simplicity to consider only the Laplace operator in
\eqref{eq.BVP}. However, all theoretical results presented here extend
to scalar, strongly elliptic, linear operators $P$ with sufficiently
regular coefficients, and even to elliptic systems, such as the system
of anisotropic elasticity \cite{MNelast}. Furthermore, we can also
treat transmission problems, for which the coefficients of $P$ are
allowed to jump across piecewise-smooth hypersurfaces, representing
interfaces, under some additional conditions \cite{BMNZ, HMN}. We
briefly discuss these extensions in Subsection \ref{ssec.extend}.

For the discretization on polyhedral domains, we build discrete spaces
$S_k \subset H^1_0(\Omega)$ and Galerkin finite element projections
$u_k \in S_k$ that approximate the solution $u$ of Equation
\eqref{eq.BVP} for $f \in H^{m-1}(\Omega)$ arbitrary. We prove that,
by using certain spaces of continuous, piecewise polynomials of degree
$m$, the sequence $S_k$ achieves {\em quasi-optimal rates of
  convergence}. More precisely we prove the existence of a constant
$C>0$, independent of $k$ and $f$, such that
\begin{equation}\label{eq.optimal.rate}
    \|u - u_k\|_{H^1(\Omega)} \le C \dim(S_k)^{-m/n}
    \|f\|_{H^{m-1}(\Omega)}, \quad u_k \in S_k,
\end{equation}
where $n=2$ or $n=3$ is the dimension of our polyhedral domain.

The contents of the paper are as follows. In the first section we
review well-posedness results in weighted Sobolev spaces on polyhedra
domains. These weighted spaces are sometimes called the \BK\ spaces.
These results are not sufficient for our applications to the Finite
Element Method in three dimensions, so in the second section we review
some additional {\em anisotropic} regularity results. These results
are used in the third section to construct a sequence of meshes that
yields $h^m$--quasi-optimal rates of convergence in three dimensions.
Finally, in the last section we discuss some of the main ingredients
that enter in the proof and which are of independent interest. These
include the Hardy-Poincar\'e inequality (which guarantees the
coercivity of our problem) and a description of the weighted Sobolev
(or \BK) spaces $\maK_{a}^m(\Omega)$, which are the natural spaces for
our well-posedness results, as the usual Sobolev spaces for a modified
metric on $\Omega$, which nevertheless is conformally equivalent to
the old one.

\vspace{.1in}

\noindent{\bf Acknowledgements:} The authors would like to thank the
organizers of the International Congress of Romanian mathematicians,
Brasov 2011, where these results were presented. A. M. would also like
to acknowledge the support of the Mathematics Department at Cornell
University, where she is currently a Michler Fellow.

%%%%%%%%%%%%%%%%%%%%%%%%%%%%%%%%%%%%%%%%%%%%%%
\section{Well posedness in isotropic weighted Sobolev spaces} 
\label{sec.one}

Using the standard notation for partial derivatives, namely $\pa_j =
\frac{\pa}{\pa x_j}$ and $\pa^\alpha = \pa_1^{\alpha_1}\ldots
\pa_n^{\alpha_n}$, for any multi-index $\alpha = (\alpha_1, \ldots,
\alpha_n) \in \ZZ_+^n$, we denote the usual Sobolev spaces on an open
set $V \subset \RR^n$ by
\begin{equation*}
    H^m(V) =\{u : V \to \CC,\ \pa^{\alpha}u \in L^2(V),\ |\alpha| \le
    m\}.
\end{equation*}
As mentioned in the introduction, the solution of our model Poisson
problem \eqref{eq.BVP} has only limited regularity in the spaces
$H^m(\Omega)$. The situation changes for the better if one considers
{\em weighted Sobolev spaces}, though. To define the weighted
analogues of these spaces, we need to introduce the notion of {\em
  singular boundary points} of the domain $\Omega \subset \RR^n$.

\subsection{Weighted Sobolev spaces}
Let $\psO \subset \pa \Omega$ be the set of singular (or non-smooth)
boundary points of $\Omega$, that is, the set of points $p \in \pa
\Omega$ such $\pa \Omega$ is not smooth in a neighborhood of
$p$. In case we consider mixed boundary conditions, the set of
  singular points includes also the set of points where the boudary
  conditions change. If, furthermore, interfaces are considered, the
  set of singular points contains the set of singular points of the
  interface, as well as the set of points where the interface touches
  the boundary. We will denote by $r_{\Omega}(x)$ the distance from a
point $x \in \Omega$ to the set $\psO$ and agree to take $r_{\Omega} =
1$ if there are no such points, i.e., if $\pa \Omega$ is smooth. For
$\mu \in \ZZ_{+}$ and $a \in \RR$, we define the weighted Sobolev
spaces
\begin{equation}\label{eq.def.wSsp0}
    \Kond{\mu}{a}(\Omega) = \{u \in L^2_{\loc}(\Omega),\,
    r_{\Omega}^{|\alpha|-a} \pa^\alpha u \in L^2(\Omega), \mbox{ for
      all } |\alpha| \le \mu\},
\end{equation}
which we endow with the induced Hilbert space norm. We note that for
$n=3$ for example and $\Omega$ a polyhedral domain in $\RR^3$, we have
that $r_{\Omega}(x) $ is the distance to the skeleton comprising the
union of the closed edges of $\partial \Omega$. Recently, general
spaces of this kind were studied by H. Amann \cite{HAmann1, HAmann2}.

Similar weighted Sobolev spaces are associated to the faces of
$\Omega$. By a {\em face}, we mean the connected components of the
boundary $\pa\Omega$ after the set of singular points is removed.  For
example for $n=3$, we define
\begin{equation*}
  \Kond{m}{a}(\pa \Omega) = \{(u_F),\,
  r_{\Omega}^{|\alpha|-a}\pa^\alpha u_F \in L^2(F)\,\},
\end{equation*}
where $|\alpha| \le m$ and $F$ ranges through the set of faces of
${\pa \Omega}$. For $s \in \RR_+$, we define the space
$\Kond{s}{a}(\pa \Omega)$ by standard interpolation.

\subsection{Well-posedness for the Poisson problem on $n$-dimensional 
polyhedral domains} The following result is proved in \cite{BMNZ}. For
simplicity, we shall assume that $\Omega$ has no cracks and that there
are no vertices that touch the boundary. (That is, we shall
  consider only domains $\Omega$ that coincide with the interior of
  their closure $\overline{\Omega}$.)

\begin{theorem}\label{theorem.main} 
Let ${\Omega \subset \RR^n}$, be a bounded, curvilinear polyhedral
domain and $m \in \ZZ_+$. Then there exists $\eta_\Omega > 0$ such
that $\tilde{\Delta}(u) = (\Delta u, u \vert_{\pa \Omega})$ defines an
isomorphism
\begin{equation*}
    \tilde{\Delta} : \Kond{m + 1}{a+1}(\Omega) \to \Kond{m - 1}{a -
      1}(\Omega) \oplus \Kond{m + 1/2}{a+1/2}(\pa \Omega),
\end{equation*}
for all $|a|< \eta_{\Omega}$. If $m = 0$, the solution $u$
corresponding to the data $(f, 0) \in \Kond{- 1}{a - 1}(\Omega) \oplus
\Kond{ 1/2}{a+1/2}(\pa \Omega)$ is also the solution of the associated
variational problem.
\end{theorem}

This theorem amounts to the well-posedness of the boundary value
\eqref{eq.BVP} on $n$-dimensional polyhedral domains. For $n=2$ (in
which case $\Omega$ is a polygonal domain), this result is due to
Kondratiev \cite{Kondratiev67}, in which case $\eta_{\Omega} =
\frac{\pi}{\alpha_{MAX}}$, where $\alpha_{MAX}$ is the measure in
radians of the maximum angle of $\Omega$. For $n=3$, this result was
proved in \cite{3D1}. For later applications, we shall need the
  following result.

\begin{theorem}\label{theorem.ext}
The results of Theorem \ref{theorem.main} remain true for infinite
angles in two dimensions, infinite polyhedral cones in three
dimensions, and infinite dihedral angles in three dimensions.
\end{theorem}

The proof of this theorem proceeds along the lines of the proof
  of Theorem \ref{theorem.main} in \cite{3D1} or \cite{BMNZ}. A first
  difference to remark is that the Hardy-Poincar\'e
  inequality 
%%CHECK
does not hold for the whole domain. 
Then the
  ``desingularization'' $\Sigma(\Omega)$ has to involve the directions
  at infinity also in the case of an angle or a cone. In case of a
  dihedral angle $D_\alpha = \{0 < \theta < \alpha\}$, in
  cylindrical coordinates $(r,\theta,z)$, one has to consider the
  also the two point compactification of the edge. In particular,
\begin{equation}
  \Sigma(D_\alpha) = [0, \alpha] \times [0, \infty] \times [-\infty,
    \infty].
\end{equation}

\subsection{Extensions} \label{ssec.extend}
Theorem \ref{theorem.main} above was extended in several ways. First,
the proof applies with almost no change if mixed boundary value
problems are considered, provided that no adjacent faces are endowed
with Neumann boundary conditions. We do allow, however, different
boundary conditions on the same face. We treat the points where the
boundary conditions change similarly to the non-smooth boundary
points, as solutions exhibit a similar singular behavior in this case.

As already mentioned in the Introduction, we can more generally
consider a general uniformly strongly elliptic differential operator
of the form 
\begin{equation}
  L u = -\sum_{ij} \pa_i(a_{ij} \pa_j u) + cu, \mbox{ with } c \ge 0.
\end{equation}
(Recall that $L$ is uniformly strongly elliptic if, and only if,
there exists $C>0$ such that $\sum_{ij} a_{ij} t_i t_j \ge C \sum_i
t_i^2$, for all $(t_i) \in \RR^n$.) 

We can also include certain transmission or interface problems.  More
precisely, we now assume that our domain $\Omega$ can be written as a
union of curvilinear polyhedral domains $\Omega_j$ with disjoint
interiors:\ $\overline{\Omega} = \cup_{j=1}^K \overline{\Omega}_j$.
Let $\Gamma := \cup_{j=1}^K \pa \Omega_j \smallsetminus \pa \Omega$ be
the interface. We assume that $\Gamma$ is smooth and assume further
that no adjacent faces of the $\Omega_j$'s are endowed with Neumann
boundary conditions. We do allow $\Gamma$ to touch the boundary of
$\Omega$.  We can then extend the result of Theorem \ref{theorem.main}
by using instead the {\em broken weighted Sobolev spaces} $\hat
\maK_{a}^{m}(\Omega)$, defined by
\begin{equation}\label{def.broken.SS}
  \hat \maK_{a}^{m}(\Omega) := \oplus_{j=1}^{K}
  \maK_{a}^{m}(\Omega_j).
\end{equation}
We observe that, if there is no interface, $\hat
\maK_{a}^{m}(\Omega)=\maK_{a}^{m}(\Omega)$.  We let $\pa_D \Omega$ be
the part of the boundary with Dirichlet boundary conditions, which we
assume to be a closed subset of the boundary, and let $\pa_N \Omega :=
\pa \Omega \smallsetminus \pa_D \Omega$. We denote the outer
  normal vector to $\Omega$, which is defined a.e. on $\pa\Omega$, by
  $\nu$, and the {\em conormal derivative} associated to the operator
  $L$ by \ $D^L_\nu = \sum_{ij} \nu_i a^{ij}(x) \pa_j$.  Let $\tilde L
  (u) = (L\, u, u\vert_{\pa_D \Omega}, D^L_{\nu} u\vert_{\pa_N
    \Omega})$. Our most general result in $n$ dimensions states that
  for $m \ge 1$ \ $\tilde L$ is an isomorphism (see \cite{BMNZ}):
\begin{equation}
    \tilde{L} :  \maD_a \to \hat \maK^{m - 1}_{a - 1}(\Omega) \oplus
    \Kond{m + 1/2}{a+1/2}(\pa_D \Omega) \oplus \Kond{m - 1/2}{a -
      1/2}(\pa_N \Omega),
\end{equation}
where 
\begin{equation} \label{eq.domaindef}
 \maD_{a} := \{ u \in \hat \maK^{m+1}_{a+1}(\Omega)\cap \maK^{1}_{a+1}
 (\Omega),\ u^+ = u^-, \ D^{L+}_\nu u = D^{L-}_\nu u \, \text{ on } \,
 \Gamma\, \},
\end{equation}
and the subscript $\pm$ refers to non-tangential limits to each side
of the interface. The conormal derivative is defined in the sense of
the trace a.e. on $\pa\Omega$.

Let us mention that the interface $\Gamma$ will separate different
faces where it touches the boundary, and hence we assume that these
faces are not both endowed with Neumann boundary conditions.

For elasticity with mixed boundary conditions, a similar result is
obtained by Mazzucato and Nistor in \cite{MNelast}. The results in
\cite{MNelast} also extend to interface problems under the same
assumptions (no adjacent faces with Neumann boundary conditions and a
smooth interface) using the methods as in \cite{BMNZ} and in
\cite{MNelast}.  More precisely, we use  Korn's inequality to
  obtain local regularity results (no weighted spaces). This
  applies, in particular, to interface problems. There the additional
  regularity is proved as for the additional regularity at the
  boundary for smooth domains. See \cite{NistorSchwab} for a proof of
  the additional regularity at the boundary for systems that extends
  to interface problems. Once one has the local regularity results,
  the global regularity results in {\em weighted} spaces is proved as
  in \cite{MNelast} using suitable partitions of unity. The
  solvability in $H^1$ is an immediate consequence of Korn's
  inequality and of the Hardy-Poincar\'e inequality. Combining
  regularity with solvability in $H^1$ yields solvability in higher
  weighted Sobolev spaces $\maK_{a+1}^{m+1}(\Omega)$.

Other regularity results go toward analytic regularity using countably
normed spaces as in the work of Babu\v{s}ka-Guo \cite{BabuGuo,
  BabuGuo2}, and Costabel, Dauge, and Nicaise \cite{CDN}. See the
Introduction for more references. It would be interesting to extend
these results to the de Rham complex \cite{ArnoldActa, ArnoldBull}.

\subsubsection{Adjacent Neumann faces and non-smooth interfaces in 2D}
The assumption that no vertex $P$ be the common point of two adjacent
faces with Neumann boundary conditions or the assumption that $\Gamma$
be smooth at any interface point $P$ are both equivalent to the fact
that the function constant equal to one not be a singular function at
that singular point $P$. This assumption is necessary, because, if it
is not satisfied, the relevant operator $\tilde L$ is not even
Fredholm for the value $a = 0$ and it is also not invertible for any
$a \in \RR$. However, this assumption is not realistic in practice
and, it turns out, not even necessary for designing graded meshes that
yield quasi-optimal rates of convergence \cite{HMN}.

To obtain a well-posedness result for interface problems in 2D, we can
proceed as follows \cite{HMN}. Let $\chi_P$ be a smooth function that
is equal to 1 near each singular point $P$ that is either a point
where we have Neumann-Neumann conditions or a non-smooth interface
point satisfying respectively Neumann or periodic boundary conditions
on the sides at $P$. This includes points $P$ that belong to more
than two of the subdomains $\overline{\Omega}_j$ (so called {\em
  multiple junction points}). We assume the $\chi$'s have disjoint
supports. Let $W_s$ be the linear span of the functions $\chi_P$.
The choice of boundary conditions or the introduction of additional
singular points to a polygonal domain define a {\em polygonal structure}
on $\Omega$, see \cite{HMN} for details.

\begin{theorem}\label{theorem.interface} 
Let $\Omega$ be a domain with a polygonal structure. The there exists
$\eta > 0$ such that, for all any $0 < a < \eta$ and $m \in \ZZ_+$,
the map 
\begin{equation*}
    \tilde{L} :  \maD_a + W_s  
    \to  \hat \maK^{m-1}_{a-1}(\Omega) \oplus
    \maK^{m + 1/2}_{a + 1/2}(\pa_D \Omega) \oplus \maK^{m - 1/2}_{a
      -1/2}(\pa_N \Omega),
\end{equation*}
with $\maD_a$ given in \eqref{eq.domaindef}, is an isomorphism.
\end{theorem}

The proof requires the calculation of the index of the operator
$\tilde L$ acting on $\hat \maK^{m+1}_{a+1}(\Omega)\cap \maK^{1}_{a+1}
(\Omega)$. Note that our result is not valid for $a = 0$. We expect a
similar result in 3D. 

Theorem \ref{theorem.interface} can be used to justify the
construction of a sequence of meshes (in 2D) that yields quasi-optimal
$h^m$ rates of convergence for transmission problems with non-smooth
interfaces (and even with multiple junctions) and problems with
adjacent Neumann-Neumann corners in 2D. See \cite{junping} for
additional issues related to the regularity and numerical methods for
interface problems. We notice that the resulting sequence of meshes is
the same for all 2D problems on polygonal domains (with or without
interfaces or Neumann-Neumann corners), although the theoretical PDE
result (or {\em a priori} estimates) are different in these two cases.

\subsection{The domain of $\Delta$ on concave polygons}
Let us mention that the method used to obtain Theorem
\ref{theorem.interface} can be used to describe the domain
$\maD(\Delta)$ of the Friedrichs extension of the Laplace operator on
$\Omega$ with zero boundary conditions. First of all, the form
associated to $\Delta$, namely $B(u, v) = (\nabla u, \nabla v)$, $u,
v$ zero on the boundary, defines the so called {\em energy
  norm}:\ $|u|_{H^1(\Omega)} = B(u, u)^{1/2}$. The completion of
$\CIc(\Omega)$ in the energy norm is $H^1_0(\Omega)$. The proofs in
\cite{BNZ1, 3D1} show that $H^1_0(\Omega) = \maK_{1}^{1}(\Omega) \cap
\{ u \vert_{\pa \Omega} = 0 \}$, with equivalent norms. The domain of
the Friedrichs extension of the Laplacian $\Delta$ is then
\begin{equation}
  \maD(\Delta) = \{u \in H^1_0(\Omega),\ \Delta u \in L^2(\Omega) \}.
\end{equation}
If $\Omega$ is convex, then it is known that $\maD(\Delta) =
H^2(\Omega) \cap H^1_0(\Omega)$. This is however not true if $\Omega$
is concave. To describe $\maD(\Omega)$ in the case when $\Omega$ is
concave, let us notice that the map
\begin{equation}
  \Delta : \maK_{2}^{2}(\Omega) \cap H^1_0(\Omega) \to L^2(\Omega)
\end{equation}
is Fredholm and its index is the number of re-entrant corners by
\cite{Kondratiev67}. Let $P$ be such a re-entrant corner with angle
$\alpha_{P} > \pi$. Also, let $(r, \theta)$ be polar coordinates at
$P$ and consider the function $\phi_{P} = r^{\pi/\alpha_{P}} \sin(\pi
\theta/\alpha_{P}) \chi_{P}$, where $\chi_{P}$ is the function
considered in Theorem \ref{theorem.interface}. Let $V_s$ be the space
of linear combinations of the functions $\phi_{P}$, with $P$ a re-entrant
corner. Then one has that
\begin{equation}
  \Delta : \maK_{2}^{2}(\Omega) \cap H^1_0(\Omega) + V_s \to
  L^2(\Omega)
\end{equation}
has index zero, is injective, and hence bijective. This proves the
following result.

\begin{theorem} The domain of the Friedrichs extension of the 
Laplace operator with zero boundary conditions on a polygon $\Omega
\subset \RR^2$ is
\begin{equation*}
  \maD(\Delta) = \maK_{2}^{2}(\Omega) \cap H^1_0(\Omega) + V_s.
\end{equation*}
\end{theorem}

A similar description is available for other types of boundary
conditions. This result immediately leads to a maximal regularity
result for the heat equation on polygonal domains.

See also \cite{GM} and \cite{GM2} for related results on Friedrichs
extensions of second order elliptic operators on manifolds with
conical points.

\section{Anisotropic weighted Sobolev spaces and regularity}\label{sec.two}
The well-posedness result of the previous section are not enough to
establish quasi-optimal rates of convergence in $3D$. We need
additional regularity along the edges, as follows.  Let $u$ be the
solution of problem \eqref{eq.BVP} with $f \in H^{m-1}(\Omega)$ and
$g=0$.  We observe that this assumption is stronger than assuming that
$f$ is in a weighted Sobolev space of the form $\maK^{m-1}_{a-1}$ for
$|a|$ small. We will need to take advantage of this additional
regularity of $f$, which leads to improved regularity for $u$
  along the edges. We encode this additional regularity by introducing
  new {\em anisotropically} weighted spaces. 

We assume first that the domain $\Omega$ is a dihedral angle with axis
along the $z$-coordinate axis, $ D_\alpha = \{0 < \theta < \alpha\}$,
using cylindrical coordinates $(r,\theta,z)$. We further assume that
$f \in H^{m-1}(D_\alpha)$. Then $f \in
  \Kond{m-1}{a-1}(D_\alpha)$, and hence
\begin{equation}\label{eq.stage1}
  u \in \Kond{m+1}{a+1}(D_\alpha)
\end{equation}
for positive and small enough $a$, by Theorem
\ref{theorem.ext}. Hence, $\pa_z u \in \Kond{m}{{a}}(D_\alpha).$
However, we also have $\Delta \pa_z u = \pa_z \Delta u = \pa_z f \in
H^{m-2}(\Omega).$ Then, using Theorem \ref{theorem.main} which extends
to this setting, we also obtain that
\begin{equation}\label{eq.stage2}
	\pa_z u \in \Kond{m}{a+1}(\Omega),
\end{equation}
a better estimate than in Equation \eqref{eq.stage1}.  These
calculations suggest that we introduce
a scale of spaces $\maD_a^m$, $m\in \ZZ_+$, as follows:
\begin{align*}
  \maD_a^1(D_\alpha) &:= \Kond{1}{1}(D_\alpha), \\ \maD_a^m(D_\alpha)
  &:=\{u \in \Kond{m}{a}(D_\alpha),\ \pa_z u \in
  \maD_a^{m-1}(D_\alpha) \}.
\end{align*}
The spaces $\maD_a^1$ are thus independent of $a$. 

We assume next that the domain $\Omega$ is a cone $\maC$ centered at
the origin. We let $\rho(x) = |x|$, the distance from $x$ to the
origin, and define
\begin{equation*}
 \maD_a^1(\maC) := {\rho^{a-1}}\Kond{1}{1}(\maC)
    =\{{\rho^{a-1}}v,\ v \in \Kond{1}{1}(\maC)\}.
\end{equation*}
To introduce the spaces $\maD_a^m(\maC)$ for $m \ge 2$, we shall need
to consider the vector field $\rho \pa_\rho := x \pa_x + y \pa_y + z
\pa_z$, which is the infinitesimal generator of dilations centered at
the vertex of the cone. We then define by induction
\begin{equation*}
 \maD_a^m(\maC) := \{u \in \Kond{m}{a}(\maC),\ {\rho \pa_ \rho} (u)
 \in \maD^{m-1}_{a}(\maC) \}, \quad m \ge 2.
\end{equation*}

For a general bounded polyhedral domain $ {\Omega }$, we define the
{\em anisotropic weighted Sobolev spaces} $\maD_a^m(\Omega)$ by
localizing around vertices and edges, using as models cones and
dihedral angles respectively, such that away from the edges these
spaces coincide with the usual Sobolev spaces $H^m$. Then, we have the
following regularity result \cite{3D2}:

\begin{theorem}\label{theorem.anisotropic} 
Let $f \in H^{m-1}(\Omega)$, with $m \ge 1$. Then there exists
$\eta_{\Omega, a} > 0$ such that the Poisson problem \eqref{eq.BVP}
with $g=0$ has a unique solution $u \in \maD^{m+1}_{a+1}(\Omega) $ for
any $0 \le a < \eta=\eta_\Omega $ and
\begin{equation*}
 \|u\|_{\maD_{a+1}^{m+1}(\Omega)} \le
    C_{\Omega, a} \|f\|_{H^{m-1}(\Omega)}.
\end{equation*}
\end{theorem}

See \cite{BKP, BCD, CDN, kellogg1} for related results.

\section{Quasi-optimal $h^m$-mesh refinement}\label{sec.three}
We describe in this section a strategy to obtain quasi-optimal
$h^m$-mesh refinement. We follow \cite{3D2}, from where the pictures
are taken. The theoretical justification of this construction is based on
the anisotropic regularity result of the previous section, Theorem
\ref{theorem.anisotropic}. 
Given a bounded polyhedral domain $\Omega$ and a parameter $\kappa
\in (0, 1/2]$, we will provide a sequence $\maT_n$ of decompositions
of $\Omega$ into finitely many tetrahedra, such that, if $S_n$ is the
finite element space of continuous, piecewise polynomials on
  $\maT_n$, 
then  is the Lagrange interpolant of $u$ of
  order $m$, $u_{I, n}$, has ``quasi-optimal'' approximability properties. The
  result can be formulated as follows:

\begin{theorem}\label{theorem.interp}\
Let $a \in (0, 1/2]$ and $0 < \kappa \le 2^{-m/a}$. Then there exists
  a sequence of meshes $\maT_n$ and a constant $C > 0$ such that, for
  the corresponding sequence of finite element spaces $S_n$, we have
\begin{equation*}
    |u - u_{I,n}|_{H^1(\Omega)} \le C 2^{-km}
    \|u\|_{\maD_{a+1}^{m+1}(\Omega)},
\end{equation*}
for any $u \in \maD_{a+1}^{m+1}(\Omega)$, $u\vert_{\pa\Omega} = 0$,
and for any $k \in \ZZ_+$.
\end{theorem}

Theorem \eqref{eq.optimal.rate} is now a direct consequence
of Theorem \ref{theorem.interp}.

\subsection{Refinement Strategy}
Our refinement strategy will first generate a sequence of
decompositions $\maT_n'$ of $\Omega$ in tetrahedra and triangular
prisms, while our meshes $\maT_n$ will be obtained by further dividing
each prism in $\maT_n'$ into three tetrahedra. We now explain how the
divisions $\maT_n'$ are defined inductively, $\maT_{n+1}'$ being a
refinement of $\maT_{n}$ in which each edge is divided into two
(possibly unequal) edges.

To define the way each edge of $\maT_n'$ is divided, we need to
establish a hierarchy of the nodes of $\maT_n'$. Therefore, given a
point $P \in \overline{\Omega}$, we shall say that $P$ is of {\em
  type} {\bf V} if it is a vertex of $\Omega$; we shall say that $P$
is of {\em type} {\bf E} if it is on an open edge of
$\Omega$. Otherwise, we shall say that it is of {\em type} {\bf S}
(that is, a ``smooth'' point). The type of a point depends only on
$\Omega$ and not on any partition or meshing. The initial \ttra\ will
consist of edges of type {\bf VE, VS, ES, EE:=E${}^2$}, and {\bf
  S${}^2$}.  We shall assume that our initial decomposition and
initial \ttra\ were defined so that no edges of type {\bf V${}^2$ :=
  VV} are present. The points of type {\bf V} will be regarded as more
singular than the points of type {\bf E}, and the points of type {\bf
  E} will be regarded as more singular than the points of type {\bf
  S}. All the resulting triangles will hence be of one of the types
{\bf VES, VSS, ESS}. Let us notice that once our initial refinement is
fine enough, the edges of our domain will be decomposed into segments
of type $VE$ and $EE$, and the segments of type $EE$ will be
containted in triangular prisms. Therefore, we can assume that there
are no triangles of type {\bf EES}.

\begin{figure}[!htb]
\begin{tabular}{cc}
{\includegraphics[width=2in]{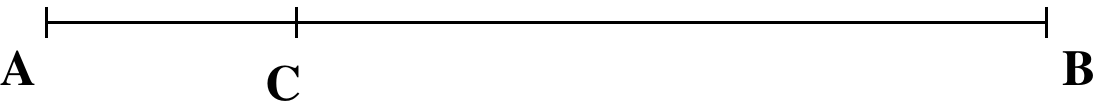}}
&
{  \includegraphics[width=2in]{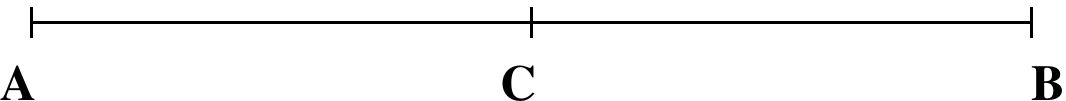} }\\
{ A more singular than B}
&
{A and B equally singular}\\
{$|AC|=\kappa |AB|, \ \kappa=1/4$}
&
{$|AC|=|AB|$}
\end{tabular}
\caption{ Edge decomposition}
\label{fig:edge}
\end{figure}

Our refinement procedure depends on the choice of a constant $\kappa
\in (0, 2^{-m/a})$, where $a> 0$ is as in Theorem
\ref{theorem.anisotropic} and $\kappa \le 1/2$. We can improve our
construction by considering different values of $\kappa$ associated to
different vertices or edges.  This generalization can easily be carry
out by the reader. See \cite{HMN} for example. Let $AB$ be a generic
edge in the decompositions $\maT_n$. Then, as part of the
$\maT_{n+1}$, this edge will be decomposed in two segments, $AC$ and
$CA$, such that $|AC|= \kappa |AB|$ if $A$ is more singular than $B$
(\ie if $AB$ is of type {\bf VE, VS}, or {\bf ES}). Except when
$\kappa = 1/2$, $C$ will be closer to the more singular point.  This
procedure is as in \cite{BNZ1}. See Figure
\ref{fig:edge}.

\begin{figure}[!htb]
\begin{tabular}{cc}
\includegraphics[width=1.8in]{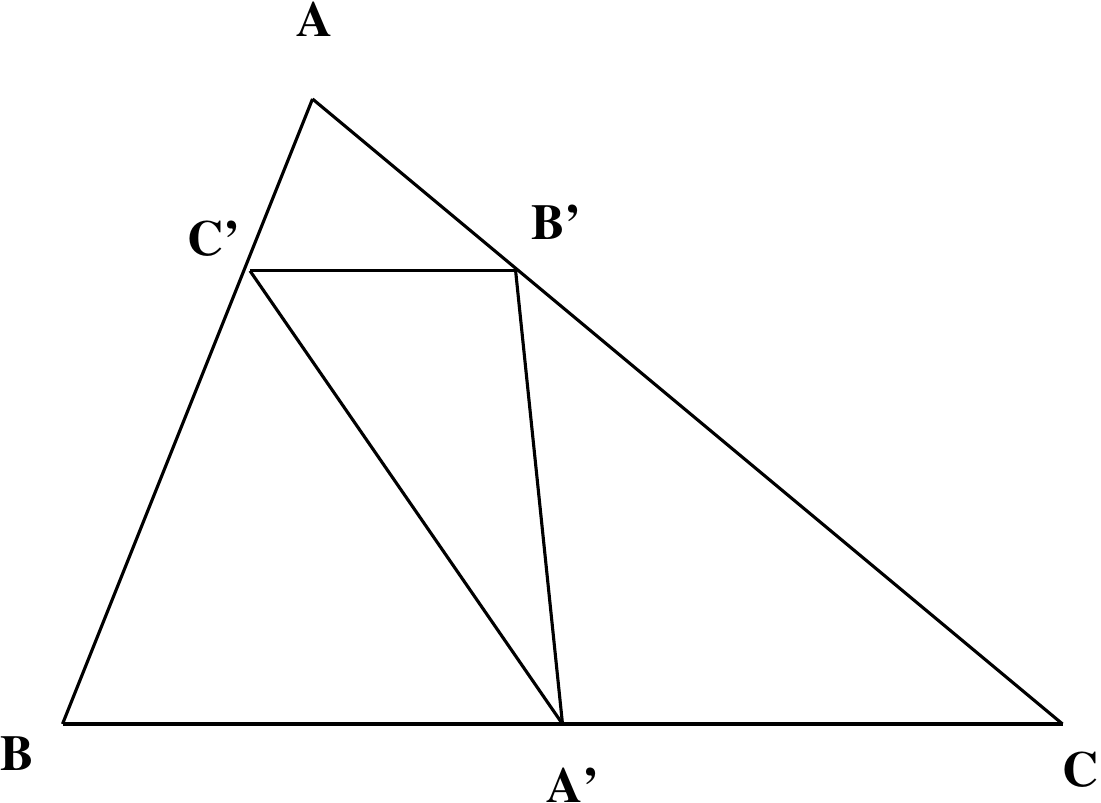}
&
\includegraphics[height=1.4in]{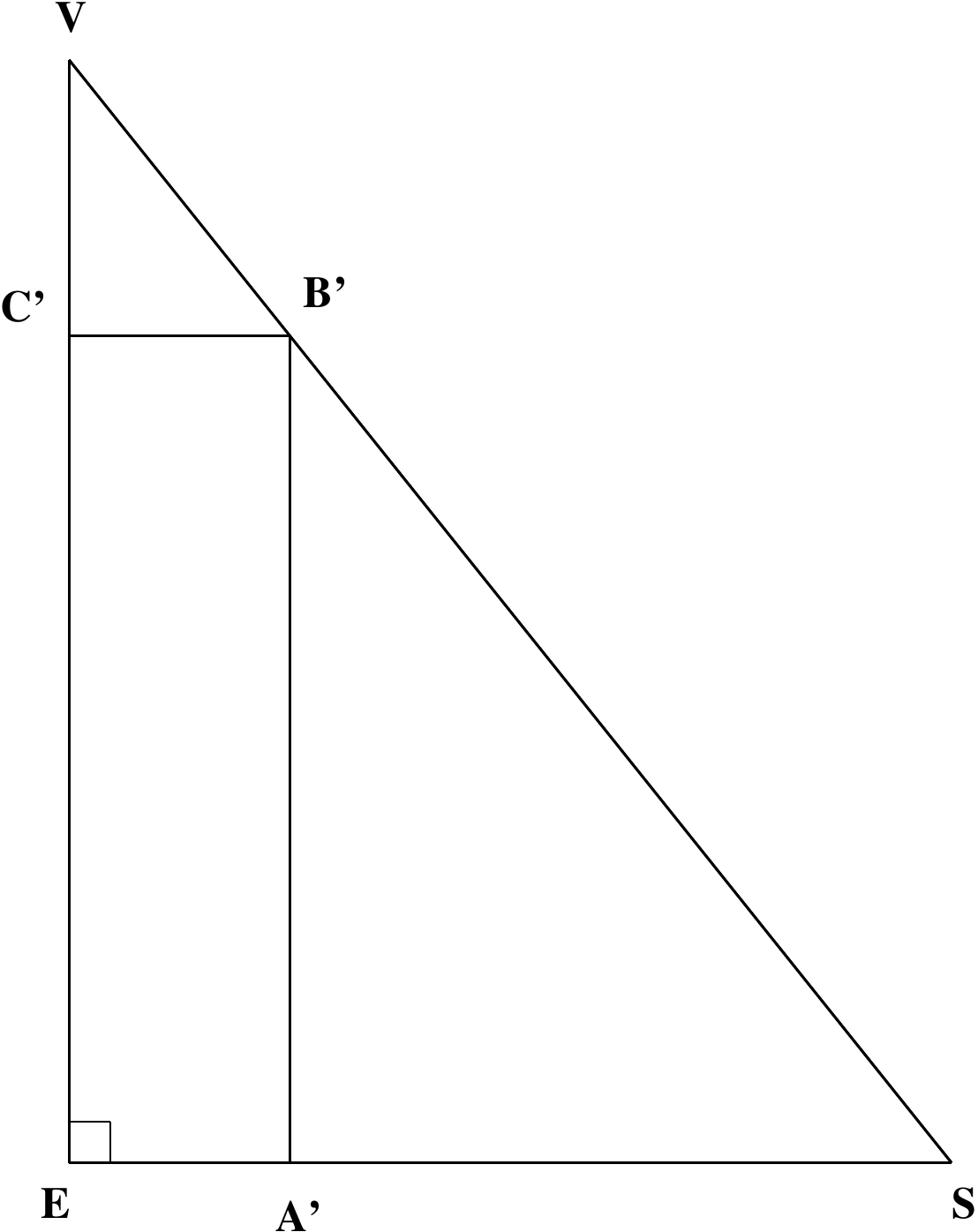}\\
{(a) $A$ of type {\bf V} or {\bf E}}
&
{VER decomposition:  $\angle E=90^o$}
\\
{(b) $B$ and $C$ of type {\bf S}, $ |A'B|=|A'C| $ }
&
{$|VC'|=\kappa |VE|,\ |VB'|=\kappa |VR|$}
\\
{ $|AC'|=\kappa |AB|,\ |AB'|=\kappa |AC|$} 
&
{$ |EA'|=\kappa |ER|$,  $A'C' $  was removed}
\end{tabular}
\caption{Triangle  decomposition, $\kappa =1/4$ }
\label{fig:VER}
\end{figure}

The above strategy to refine the edges induces a natural strategy for
refining the triangular faces. If $ABC$ is a triangle in the
decomposition $\maT_n'$, then in $\maT_{n+1}'$, the triangle $ABC$
will be divided into four other triangles, according with the edge
strategy. The decomposition of triangles of type {\bf S${}^3$} is
obtained for $\kappa = 1/2$. The type\ {\bf VSS} triangle
decomposition is described in Figure \ref{fig:VER} (a). In the
case when $ABC$ is of type {\bf VES}, however, we shall use a
different construction. Namely, in this case we remove the newly
introduced segment that is opposite to $B$, see Figure \ref{fig:VER}
(b), and divide $ABC$ into two triangles and a
quadrilateral. The resulting quadrilateral will belong to a prism in
$\maT_{n+1}'$.

\subsection{Divisions  in tetrahedra and prisms}
We now describe the construction of the sequence of the decompositions
$\maT_n'$ for $n \ge 0$. The required sequence of meshes $\maT_n$ will
be defined by dividing all the prisms in $\maT_n'$ into
tetrahedra. For the first level of semi-uniform refinement of a prism,
more details are presented in \cite{3D2}.

\begin{figure}[!htb] 
\label{fig:TWO}
\begin{tabular}{cc}
\includegraphics[width=2.2in]{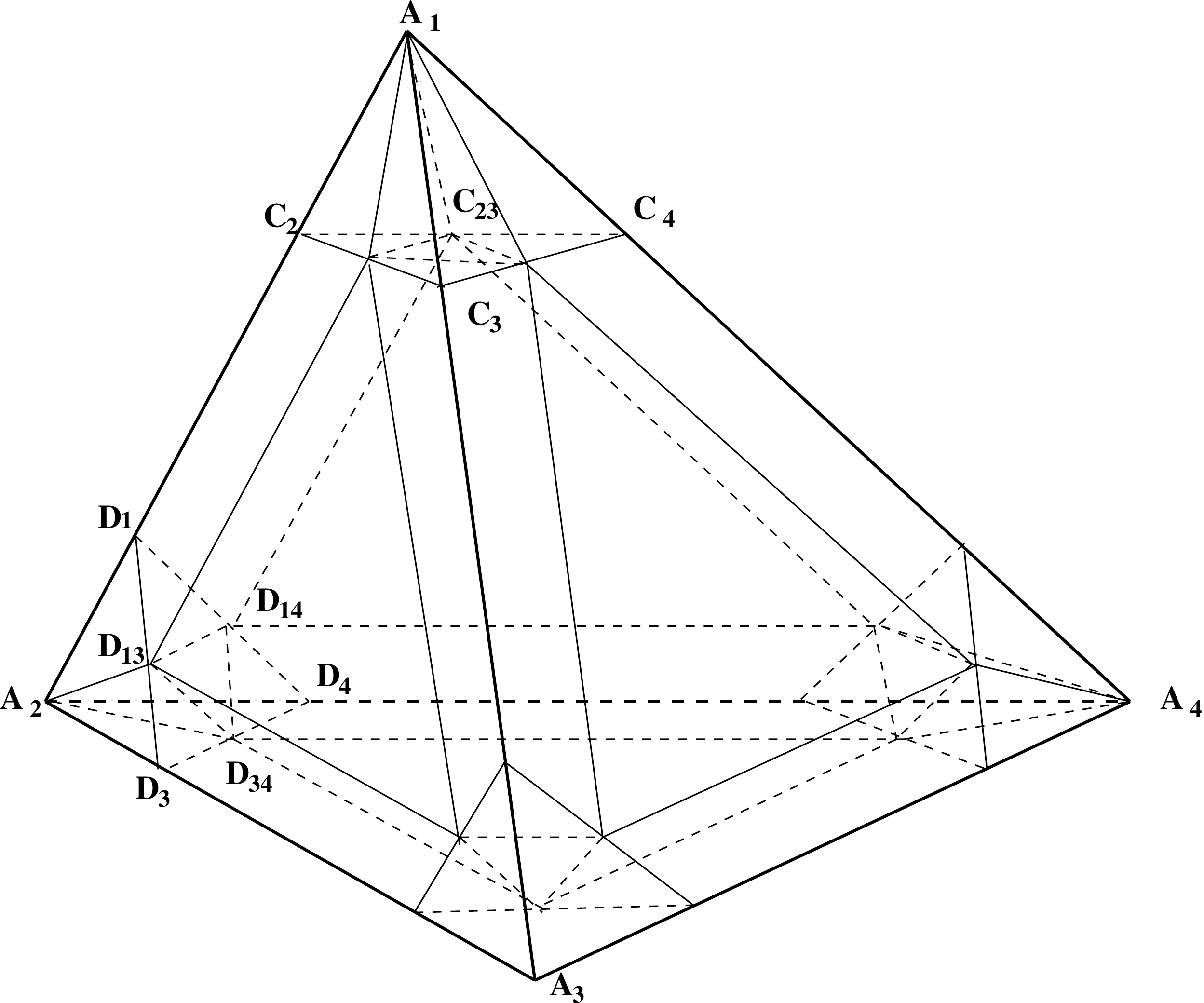} &
\includegraphics[width=1.8in]{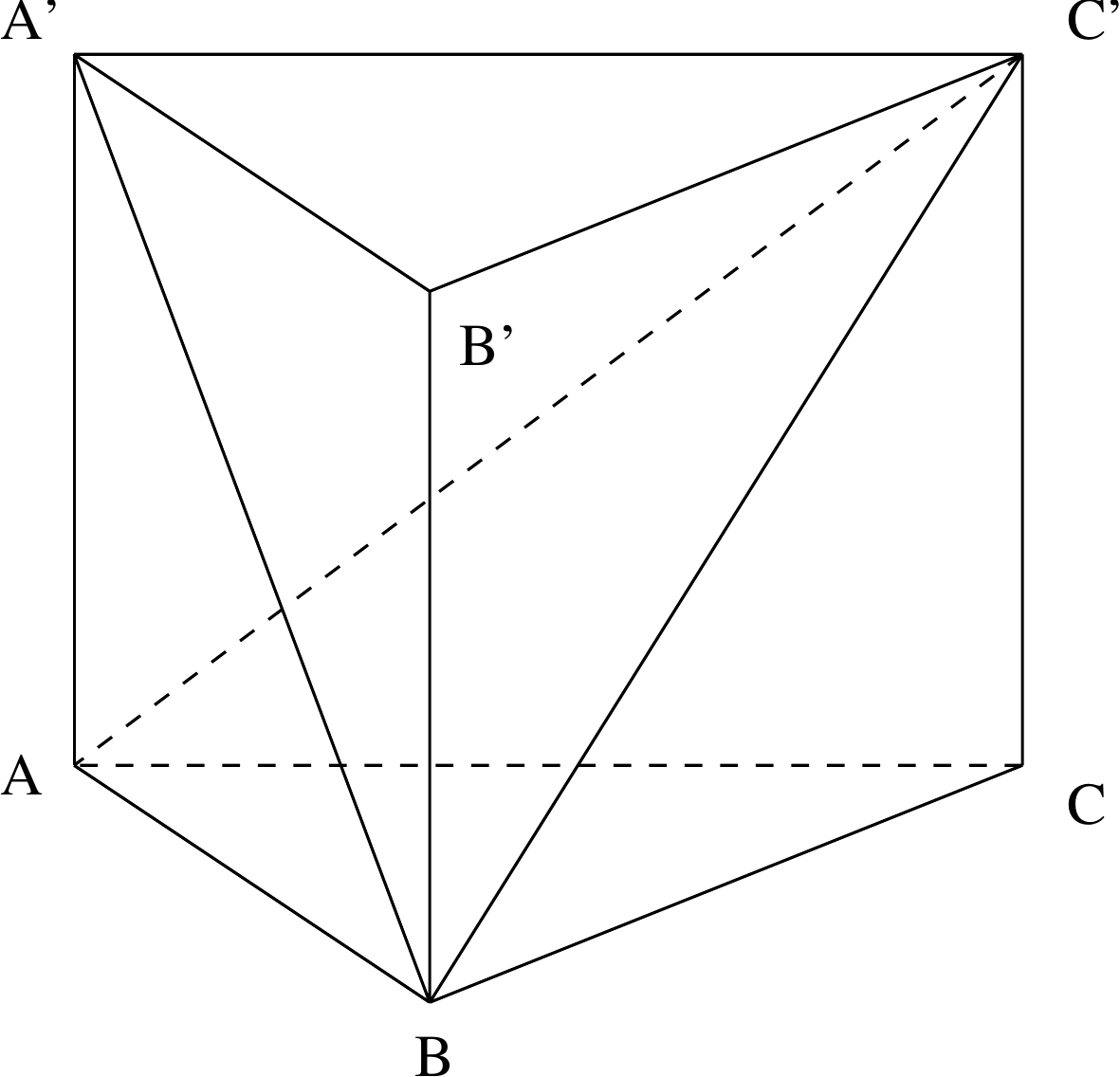}\\
(a) Initial decomposition. &
(b) Marking  a prism:
$BC' = mark$,\\
& $AA'\  ||\  BB'\  || \ CC' \perp ABC$ and
$A'B'C' $
\end{tabular}
\caption{The initial decomposition $\maT_0'$ of $\Omega$.}
\end{figure}

We start with an initial division $\maT_0'$ of $\Omega$ in straight
triangular prisms and tetrahedra of types {\bf VESS} and {\bf
  VS${}^3$}, having a vertex in common with ${\Omega }$, and an
interior region $\Lambda_0$. See Figure \ref{fig:TWO}
(a), where we have
assumed that our domain $\Omega$ is a tetrahedron. For each of the
prisms we choose a diagonal (called mark) which will be used to
uniquely define a partition of the triangular prism into
tetrahedra. We then divide the interior region $\Lambda_0$ into
tetrahedra that will match the marks.  Also, we assume that the marks
on adjacent prisms are compatible, so that the resulting meshes are
conforming. We can further assume that some of the edge points (as in Figure
\ref{fig:TWO} (b))  have been moved along the edges so that the prisms
become straight triangular prisms {\em i.e.}, the edges are
perpendicular to the bases.

The decompositions $\maT_n'$ are then obtained by induction following
the Steps {\bf 1} through {\bf 3} explained next. We assume that the
decomposition $\maT_n'$ was defined and we proceed to define the
decomposition $\maT_{n+1}'$.
\medskip

\begin{figure}[!htb] 
\label{fig:THREE}
\begin{tabular}{cc}
\includegraphics[width=2.2in]{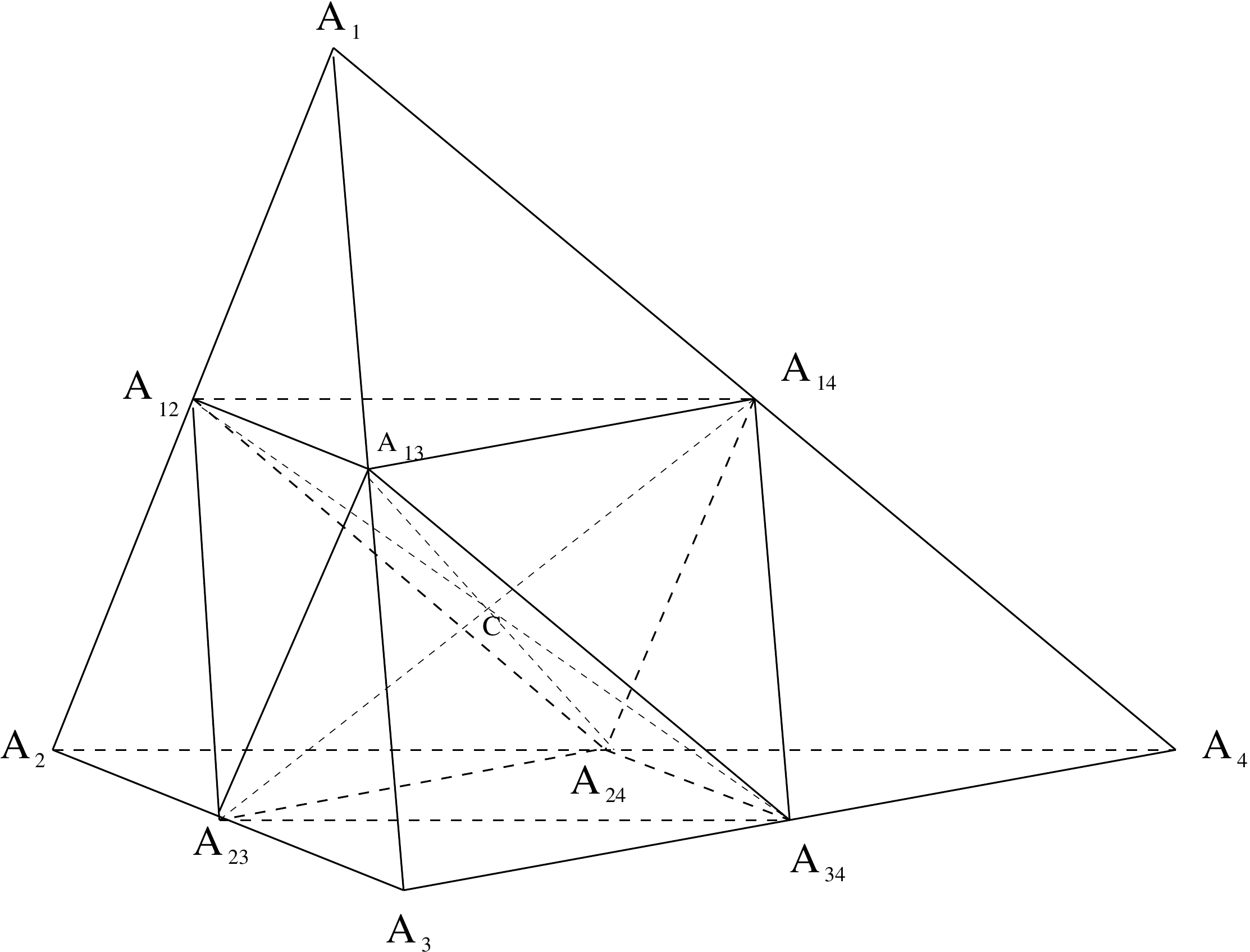} &
\includegraphics[width=1.3in]{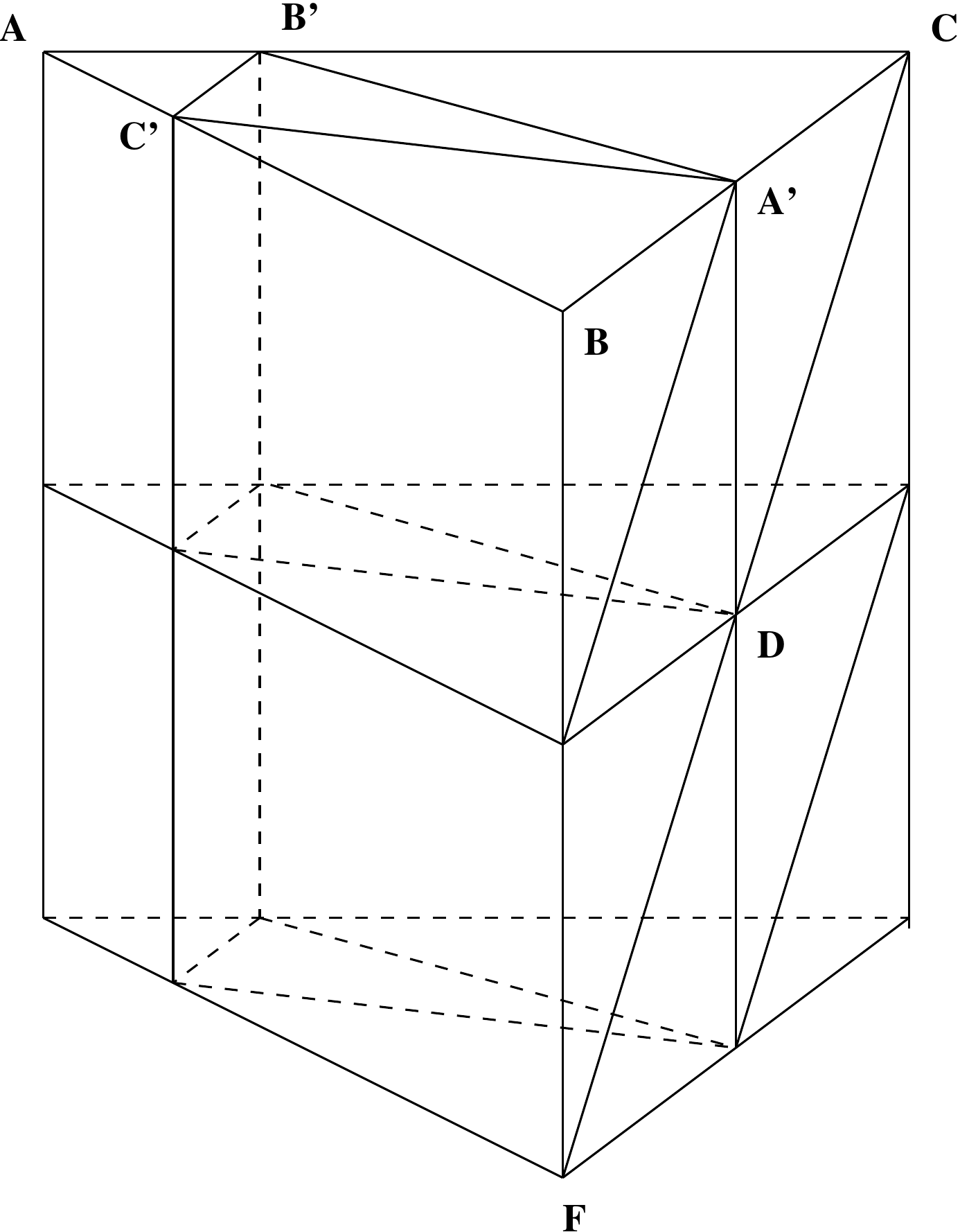} \\
(a) First level of uniform refinement &
(b) First level of semi-uniform refinement\\ 
& of a prism, $CD = mark$
\end{tabular}
\caption{First refinement $\maT_1'$.}
\end{figure}

\noindent {\bf Step 1.} The tetrahedra of type {\bf $S^4$} are refined
uniformly by dividing along the planes given by \ ${x_i + x_j =
  k/2^n}, {1 \le k \le 2^n}$, where ${x_j}$ are affine barycentric
coordinates. This refinement is compatible with the already defined
refinement procedure for the faces. See Figure \ref{fig:THREE} (a) for
$n=1$.
\smallskip

\noindent {\bf Step 2.} We perform semi-uniform refinement for prisms
in our decomposition $\maT_{n}'$ (all these prisms will have an edge
in common with $\Omega$). This procedure is shown in Figure
\ref{fig:THREE} (b).

\begin{figure}[!htb] 
\label{fig:FOUR}
\begin{tabular}{cc}
\includegraphics[width=1.8in]{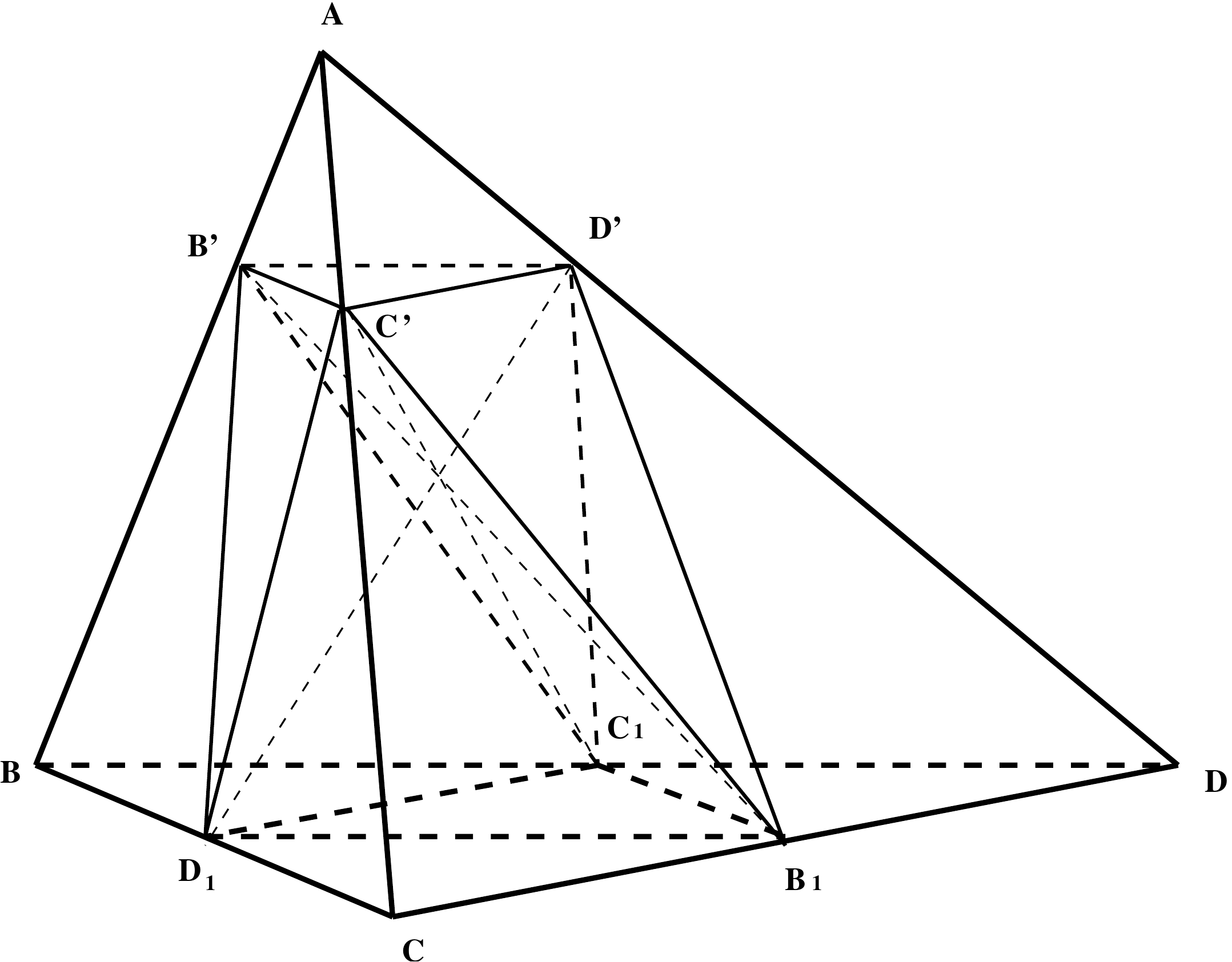} &
\includegraphics[width=1.8in]{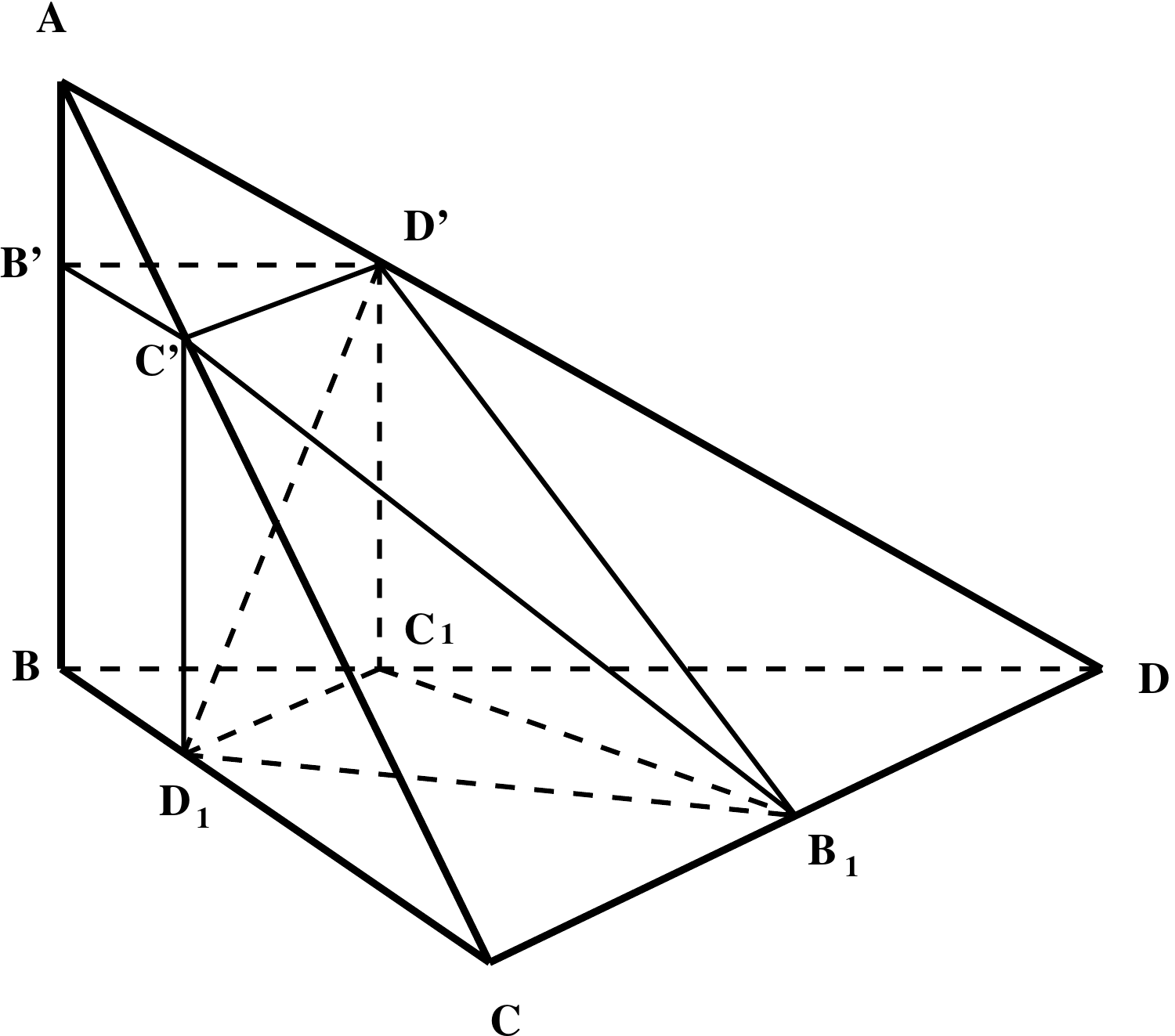} \\
(a) Vertex A of type {\bf V}, &
(b) Vertex A of type {\bf V}, B of type {\bf E},\\
B, C, D of type {\bf S} &
C, D of type {\bf S} and $D_1D'=$ mark \\
& for the prism $BD_1 C_1 D' C_1B'$
\end{tabular}
\caption{Refinement of tetrahedra of type {\bf VS${}^3$} and {\bf
    VESS}.}
\end{figure}

\noindent {\bf Step 3.}  We perform non-uniform refinement for the
tetrahedra of type {\bf VS${}^3$} and {\bf VESS}. More precisely, we
divide a tetrahedron of type {\bf VS${}^3$} into 12 tetrahedra as in
the uniform strategy, with the edges through the vertex of type {\bf
  V} divided in the ratio ${\kappa }$. We thus obtain one tetrahedron
of type {\bf VS${}^3$} and 11 tetrahedra of type {\bf S${}^4$}. (At
the next step, which yields $\maT_{n+2}'$ we iterate this procedure
for the small tetrahedron of type {\bf VS${}^3$}, while the tetrahedra
of type {\bf S${}^4$} are divided uniformly.) See Figure
\ref{fig:FOUR} (a). On the other hand, a tetrahedron of type {\bf
  VESS} will be divided it into 6 tetrahedra of type {\bf S${}^4$},
one tetrahedron of type {\bf VS${}^3$}, and a triangular prism. The
vertex of type {\bf E} of will belong only to the prism.  See Figure
\ref{fig:FOUR} (b). This refinement is compatible with the earlier
refinement of the faces.

The description of our refinement procedure is now complete.
\medskip

\subsection{Intrinsic local refinement} 
We see that one of the main features of our refinement is that each
edge, each triangle, and each quadrilateral that appears in a
tetrahedron or prism in the decomposition $\maT_n'$ is divided in the
decomposition $\maT_{n+1}'$ in an intrinsic way that depends only on
the type of the vertices of that edge, triangle, or quadrilateral. In
particular, the way that a face in $\maT_n'$ is divided to yield
$\maT_{n+1}'$ does not depend on the type of the other vertices of the
tetrahedron or prism to which it belongs. This ensures that 
\ttra\ $\maT_{n+1}$, which is obtained from $\maT_{n+1}'$ by dividing
each prism in three tetrahedra, is a conforming mesh.

% \section{Conclusion}

\section{Hardy-Poincar\'e inequality and regularity: a glimpse at the proofs}
\label{sec.four}

There are two main ingredients for the proofs of the well-posedness
results stated in the first section. One is the Hardy-Poincar\'e
inequality, which yields solvability (more precisely well-posedness)
in the $H^1$-type spaces and the second one is a regularity result,
which allows us then to obtain well-posedness in higher regularity
spaces. A third, more technical ingredient, is to describe the trace
spaces at the boundary. For this, we use the same ideas as the ones
used in the proof of regularity. We now discuss these ingredients.

\subsection{The Hardy-Poincar\'e inequality} 
Let us denote by $r_{\Omega}(x)$ the distance from $x$ to the set of
singular points in the boundary of $\Omega$.  Recall that these
singular points consist not just of the edge points, but also of the
points where the boundary conditions change and the points where the
interface touches the boundary. The following inequality is then
proved by induction \cite{BMNZ} (see \cite{3D1} for the three
dimensional case, the two dimensional case was well known, see
\cite{NP} for example).

\begin{proposition} \label{prop.HP}
Let $\Omega$ be a polyhedral domain in $\RR^n$. We assume that either
$\Omega$ is bounded, or that it is a cone or a dihedral angle. Let us
assume that the Neumann part of the boundary $\pa_N \Omega := \pa
\Omega \smallsetminus \pa_D \Omega$ contains no adjacent faces of
$\Omega$. Then there exists a constant $C_{\Omega} > 0$, which depends
only on $\Omega$ and the choice of boundary conditions such that the
following {\em Hardy-Poincare\'e} inequality holds:
\begin{equation*}
  \int_{\Omega} \frac{|u|^2}{r_\Omega^2} dx \le C_{\Omega}
  \int_{\Omega} |\nabla u|^2 dx
\end{equation*}
for any function $u \in H^1(\Omega)$ that is zero on $\pa_D \Omega$.
\end{proposition}

Let us assume that $\Omega$ is bounded. A simple consequence of the
Hardy-Poincar\'e inequality of Proposition \ref{prop.HP} is that the
spaces
\begin{equation}
  H_D^1(\Omega) = \{ u\in H^1(\Omega), u = 0 \mbox{ on } \pa_D \Omega
  \}
\end{equation} 
and 
\begin{equation}
  \maK_{1}^{1}(\Omega) \cap \{ u\in H_{loc}^1(\Omega), u = 0 \mbox{ on
  } \pa_D \Omega \}
\end{equation}
are the same and their respective norms are equivalent. Neither this
result nor the Hardy-Poincar\'e inequality are true if there exist two
adjacent faces with Neumann boundary conditions. This is the reason we
needed a different approach in Section~\ref{sec.one}.

\subsection{Sobolev spaces and regularity} 
Our definition of weighted Sobolev spaces, Equation
\eqref{eq.def.wSsp0}, is elementary. However, for the purpose of
establishing the needed properties of these spaces, it is convenient
to identify them with the usual Sobolev spaces associated to a
different metric on $\Omega$.

To this end, let us recall from \cite{BMNZ} that a {\em stratified
  curvilinear polyhedral domain} $\Omega$ is an open subset of a
Riemannian manifold $(M, g)$ of dimension $d$ together with a
stratification of 
\begin{equation}
  \overline{\Omega} = \Omega^{(d)} \supset \Omega^{(d-1)} \supset
  \ldots \supset \Omega^{(1)} \supset \Omega^{(0)}.
\end{equation}
We then define stratified curvilinear polyhedral domains by induction
as follows. For $d = 0$, $\Omega$ is just a finite set of points. For
$d = 1$, $\Omega$ is a finite set of intervals. The stratum $S_0$ for
$d = 1$ will contain all the boundary points of the intervals, but may
contain also other points. For $d > 1$, we require our domain $\Omega$
to satisfy the following conditions: for every point $p \in \pa
\Omega$, there exist a neighborhood $V_p \subset M$ such that, if
$p\in \Omega^{(l)}\setminus \Omega^{(l-1)}$, $l=1,\dots, d-1$,
then there is a stratified curvilinear polyhedral domain $\omega_p
\subset S^{d-l-1}$, $\overline{\omega_p} \neq S^{d-l-1}$, and a
diffeomorphism $\phi_p : V_p \to B^{d-l} \times B^{l}$ such that
$\,\phi_p(p)=0$ and
\begin{equation}\label{eq.cond.poly}
    \phi_p(\Omega \cap V_p) = \{r x',\, 0 < r < 1,\, x' \in \omega_p\}
    \times B^{l},
\end{equation}
inducing a homeomorphism $\overline{\Omega} \cap V_p \to \{r x',\, 0
\le r < 1,\, x' \in \overline{\omega_p}\} \times B^{l}$ of stratified
spaces that is a diffeomorphism on each stratum. 

The set of singular points of $\Omega$ then consists of
$\Omega^{(n-2)}$ and is given as part of the definition of $\Omega$,
but it must contain all the geometric, intrinsic singular points of
$\pa \Omega$.  Althought we shall not need this definition here, let
us mention nevertheless that the {\em desingularization} of $\Omega$,
denoted $\Sigma(\Omega)$, is obtained by gluing in a natural way all
the sets $[0, 1) \times \overline{\omega_p} \times B^{l}$ as in
  Equation \eqref{eq.cond.poly}. The resulting set $\Sigma(\Omega)$ is
  then a manifold with corners that has a natural structure of a Lie
  manifold with boundary, in the sense of \cite{AIN}. Then
  $\Sigma(\Omega) \to \Omega$ is a differentiable map that is a
  diffeomorphism outside the set of singular points, in
  $\Sigma(\Omega)$ the set of singular points being the set of points
  belonging to a face of codimension at least two.

Let $\tilde r_0(x) \ge 0$ be the distance from $x$ to the set
$\Omega^{(0)}$ if $x$. In general, the function $\tilde r_0$ will not
be smooth, we therefore replace $\tilde r_0$ with an equivalent
function $r_0$ that is smooth outside $\Omega^{(0)}$.  Therefore, we
also have that $r_0(x) > 0$, for $x \notin \Omega^{(0)}$, and that
$\tilde r_0/r_0$ and $r_0/\tilde r_0$ are bounded functions. We shall
say that $r_0$ is the {\em smoothed distance} to $\Omega^{(0)}$. We
replace then the metric $g =: g_0$ with $g_1 := r_0^{-2}g$. We repeat
this construction for the remaining non-empty strata in the increasing
order of the dimension of the strata, each time measuring distances in
the new metric. Thus $r_k$ is the smoothed distance to $\Omega^{(k)}$
in the metric $g_k$, and we let $g_{k+1} := r_{k}^{-2}g_{k}$, $k \le d
- 2$. One can prove that $g_{d-1}$ is a compatible metric on the
desingularization $\Sigma(\Omega)$ \cite{aln1, BMNZ} and hence we can
use the results on Sobolev spaces from those papers. Let $\rho := r_0
r_1 \ldots r_{d-2}$. Let us denote by $\Gamma(\overline{\Omega}, TM)$
the space of restrictions to $\overline{\Omega}$ of smooth vector
fields on $M$. The resulting structural Lie algebra of vector fields
on $\Sigma(\Omega)$ is simply $\maV = \CI(\Sigma(\Omega)) \rho
\Gamma(\overline{\Omega}, TM)$. Thus a basis of $\maV$ over
$\CI(\Sigma(\Omega))$ is given by $\{\rho \pa_i\}$. The resulting
Sobolev spaces are therefore
\begin{equation*}
  \maK_{a}^{m} (\Omega) := \{ u, \rho^{|\alpha| - a} \pa^\alpha u \in
  L^2(\Omega),\ |\alpha| \le m \} 
  = \rho^{a - n/2} H^m(\Omega, g_{d-1}),
\end{equation*}
where the space $H^m(\Omega, h)$ is the Sobolev space associated to
the metric $h$. Let $r_{\Omega}(x)$ denote the distance from $x$ to
$\Omega^{(d-2)}$. One can prove by induction that $r_{\Omega}/\rho$
and $\rho/r_{\Omega}$ are both bounded, so in the above definition of
Sobolev spaces we can replace $\rho$ with $r_{\Omega}$. See
\cite{BMNZ} for details.

The fact that the Sobolev spaces $\maK_{a}^{m}$ are associated to a
Lie manifold guarantees that Laplacian $\Delta$ satisfies elliptic
regularity in the scale of spaces $\maK_{a}^{m}(\Omega)$. To this end,
one also needs to establish that $\rho^2 \Delta - \Delta_{g_{d-1}}$ is
a lower order differential operator generated by $\maV$ and
$\CI(\Sigma(\Omega))$. We also obtain as a byproduct the fact that the
traces at the boundary of the spaces $\maK_{a}^{m}(\Omega)$ can also
be described in terms of the Sobolev spaces on $\pa \Omega$ associate
to the conformally equivalent metric $h$.

The Hardy-Poincar\'e inequality can also be interpreted in the setting
of the desigularized metric. Indeed, we have that there exists $C > 0$
such that every point of $x$ is at a distance $\le C$ to the {\em
  Dirichlet part} of the boundary of $\Sigma(\Omega)$ if, and only if,
there exist no two adjacent faces with Neumann boundary
conditions. Then, once we know that every point is at a distance $\le
C$ to the Dirichlet boundary, we can prove the Hardy-Poincar\'{e}
inequality in the usual way.

% \bibliographystyle{plain}
% \bibliography{bonp}

\def\cprime{$'$} \def\cprime{$'$} \def\cprime{$'$} \def\cprime{$'$}
  \def\cprime{$'$} \def\ocirc#1{\ifmmode\setbox0=\hbox{$#1$}\dimen0=\ht0
  \advance\dimen0 by1pt\rlap{\hbox to\wd0{\hss\raise\dimen0
  \hbox{\hskip.2em$\scriptscriptstyle\circ$}\hss}}#1\else {\accent"17 #1}\fi}
  \def\cprime{$'$} \def\udot#1{\ifmmode\oalign{$#1$\crcr\hidewidth.\hidewidth
  }\else\oalign{#1\crcr\hidewidth.\hidewidth}\fi}

\end{document}